\documentclass{amsart}
\title[Simple proofs of classical reciprocity laws]{Simple
proofs of classical explicit reciprocity laws on curves using
determinant groupoids over an artinian local ring}
\author[G.\ W.\ Anderson and F.\ Pablos Romo]{$\;$}
\thanks{The second-named
author was partially supported by DGESYC research contract
BFM2000-1327 and Castilla y Le\'on regional government contract
SA064/01.\\
2000 {\em Mathematics Subject Classification}: 14H05, 11R37,
19F15}
\date{This paper has been accepted for publication in Communications
in Algebra. It is to be hoped that it appears in 2003 or 2004.}
\newcommand{\Witt}{{\mathbb W}}
\newcommand{\PP}{{\mathcal P}}
\newcommand{\OO}{{\mathcal O}}
\newcommand{\ZZ}{{\mathbb Z}}
\newcommand{\iso}{\xrightarrow{\sim}}
\DeclareMathOperator{\Res}{{\mathrm{Res}}}

\DeclareMathOperator{\rank}{{\mathrm{rk}}}
\DeclareMathOperator{\Rank}{{\mathrm{Rk}}}
\DeclareMathOperator{\Hom}{{\mathrm{Hom}}}
\DeclareMathOperator{\ind}{{\mathrm{ind}}}
\DeclareMathOperator{\Det}{{\mathrm{Det}}}
\newcommand{\QQ}{{\mathbb{Q}}}
\newcommand{\xbold}{{\mathbf{x}}}
\newcommand{\ybold}{{\mathbf{y}}}
\newcommand{\Abold}{{\mathbf{A}}}
\newcommand{\Mbold}{{\mathbf{M}}}
\newtheorem{Lemma}[subsubsection]{Lemma}
\newtheorem{Proposition}[subsubsection]{Proposition}
\newtheorem{Theorem}[subsubsection]{Theorem}
\DeclareMathOperator{\Mat}{{\mathrm{Mat}}}
\begin{document}
\maketitle

{\center Greg W.\ Anderson\\
         School of Mathematics\\
         University of Minnesota\\
         Minneapolis, MN 55455, USA\\
         {\tt gwanders@math.umn.edu}\\
         and\\
         Fernando Pablos Romo\\
         Departamento de Matem\'{a}ticas\\
         Universidad de Salamanca\\
         Plaza de la Merced 1-4\\
         37008 Salamanca, Espa\~{n}a\\
         {\tt fpablos@usal.es}\\
         \;}

\begin{quote}{\small {\bf Abstract.}
The notion of determinant groupoid is a natural
outgrowth of the theory of the Sato Grassmannian and
thus well-known in mathematical physics. We briefly sketch
here a version of the theory of determinant groupoids over an artinian
local ring, taking pains to put the theory in  a simple concrete form
suited to number-theoretical applications. We
then use the theory to give a simple proof of a reciprocity law for the
Contou-Carr\`{e}re symbol. Finally, we explain how from the latter to
recover various classical explicit reciprocity laws on nonsingular complete
curves over an algebraically closed field, namely
sum-of-residues-equals-zero, Weil reciprocity,  and an explicit reciprocity
law due to Witt. Needless to say, we have been much influenced by the work
of Tate  on sum-of-residues-equals-zero and   the work of
Arbarello-DeConcini-Kac on Weil reciprocity. We also build in an essential
way on a previous work of the second-named author.}
\end{quote}
\bigskip

\section{Introduction}

    In 1968 J. Tate \cite{Tate} gave a definition
of the residues of differentials on a curve in
terms of traces of certain linear operators on infinite-dimensional
vector spaces. Further, Tate deduced the residue theorem
(``sum-of-residues-equals-zero'') on a nonsingular complete curve $X$
 from the finite-dimensionality of the
cohomology groups
$H^{0}(X,\OO_X)$ and $H^{1}(X,\OO_X)$. This work of Tate has been
enormously influential.

   In 1989 E.\ Arbarello, C.\ De Concini and
V.\ G.\ Kac \cite{AdCK} interpreted the tame symbol  at a point of 
a complete nonsingular
algebraic curve
$X$ over an algebraically closed field as a commutator in a
certain central extension of groups and then, in the style of Tate,
deduced a reciprocity law  on
$X$ for the tame symbol (``Weil reciprocity'')
from the finite-dimensionality of
$H^0(X,\OO_X)$ and $H^1(X,\OO_X)$.
Recently the second-named author of this paper has provided an
interpretation  \cite{PablosRomo} of the central
extension of \cite{AdCK} in terms of determinants associated to
infinite-dimensional vector subspaces
 valid for curves over a perfect
field. The logical organization of this paper to a significant extent
parallels that of \cite{PablosRomo}.

    In 1994 C. Contou-Carr\`{e}re \cite{ContouCarrere} defined a
natural transformation greatly generalizing the tame symbol. In the
case of an artinian local base ring $k$ with 
maximal ideal $m$,
the natural transformation takes the following form. 
Let 
$f,g\in k((t))^\times$ be given, where $t$ is a variable.
(Here and below
$A^\times$ denotes the multiplicative group of a ring $A$ with unit.)
It is possible in exactly one way to write 
$$\begin{aligned} f &=
a_0\cdot t^{w(f)} \cdot {\prod}_{i=1}^{\infty} (1-a_{i}t^{i}) \cdot
{\prod}_{i=1}^{\infty} (1-a_{-i}t^{-i}) \\ g &= b_0\cdot t^{w(g)}
\cdot {\prod}_{i=1}^{\infty} (1-b_{i}t^{i}) \cdot
{\prod}_{i=1}^{\infty} (1-b_{-i}t^{-i})\end{aligned}$$\noindent with
$w(f),w(g)\in \ZZ$, 
$a_{ i},b_{i}\in k$ for 
$i>0$, $a_0, b_0 \in k^\times$, $a_{-i}, b_{-i}\in m$ for $i>0$,
and $a_{-i}=b_{-i}=0$ for $i\gg 0$.
By definition the value of the {\em Contou-Carr\`{e}re
symbol} is
$$\langle f,g\rangle  := (-1)^{w(f)w(g)}\frac
{a_0^{w(g)}{\prod}_{i=1}^{\infty}{\prod}_{j=1}^{\infty}\big (1 -
a_i^{j/(i,j)}b_{-j}^{i/(i,j)}\big
)^{(i,j)}}{b_0^{w(f)}{\prod}_{i=1}^{\infty}{\prod}_{j=1}^{\infty}\big
(1 - a_{-i}^{j/(i,j)}b_{j}^{i/(i,j)}\big )^{(i,j)}} \in
k^\times.$$
The definition makes sense because only finitely many of the terms
appearing in the infinite products differ from 1. The symbol
$\langle\cdot,\cdot\rangle$ is clearly antisymmetric and, although it
is not immediately obvious  from the definition, also
bimultiplicative.

If $k$ is a field, and hence $m=0$, then the infinite products go
away and the Contou-Carr\`{e}re
symbol reduces to the tame symbol.
If $k=k_0[\epsilon]/(\epsilon^3)$, where $k_0$ is a field,
then
$$\langle 1-\epsilon f,1-\epsilon g\rangle\equiv 1
-\epsilon^2\Res_{t=0}( g\,df)\bmod{\epsilon^3}$$
for all $f,g\in k_0((t))$,
and so the Contou-Carr\`{e}re symbol also contains the residue as a
special case. 
If $k$ is a
$\QQ$-algebra and $f\in 1+m((t))$, then
$$\langle f,g\rangle=\exp(\Res_{t=0}\log f\cdot d\log g).$$
This last formula renders the bimultiplicativity of
the Contou-Carr\`{e}re symbol at least plausible and motivates the
definition.

    The main aims of this paper are (i) to interpret the Contou-Carr\`{e}re
symbol $\langle f,g\rangle$---up to signs---as a commutator of liftings 
of $f$ and $g$
to a certain central extension of a group containing $k((t))^\times$
(see Thm.~\ref{Theorem:CommutatorInterpretation})
and then (ii) to exploit the commutator interpretation to prove in the
style of Tate a reciprocity law for the Contou-Carr\`{e}re symbol
on a nonsingular complete curve defined over an algebraically closed field
(see Thm.~\ref{Theorem:CCR}).
The commutator interpretation of the Contou-Carr\`{e}re 
symbol provided here formally
resembles the commutator formula \cite[Prop.\ 3.6]{SegalWilson}
stated in Segal-Wilson.

In more detail, the general reciprocity law proved here takes the
following form. Let $F$ be an algebraically closed field.
Let $X/F$ be a complete nonsingular curve.
Let $S$ be a finite nonempty set of (closed) points of $X$.
For any ring or
group
$A$, put
$A^S:=\{(a_s)_{s\in S}\vert a_s\in A\}$.  By choosing uniformizers
at each point belonging to
$S$, identify
$R_0:=H^0(X\setminus S,\OO_X)$ with an $F$-subalgebra 
of $F((t))^S$.
Further, suppose  now that the artinian local ring
$k$ considered above is a finite $F$-algebra. Put
$R:=R_0\otimes_Fk$ and make the evident identification of $R^\times$
with a subgroup of $k((t))^{\times S}$. We prove
that
$$\prod_{s\in S}\langle f_s,g_s\rangle=1$$
for all $f,g\in R^\times$.  With $k=F$ we get back Weil
reciprocity. With
$k=F[\epsilon]/(\epsilon^3)$ we get back sum-of-residues-equals-zero.
With $F$ of characteristic $p>0$ and
$k=F[\epsilon]/(\epsilon^{p^{n-1}+1})$, we get back  an explicit
reciprocity law due to Witt \cite{Witt}. 

The general reciprocity law proved here seems in principle to be known
in the mathematical physics world, albeit no reference convenient for a
number theorist can be cited.  We claim novelty only for the
simplicity and directness of our approach to the reciprocity law.
We hope to popularize the reciprocity law among
number theorists, and expect it to have many applications.

The backdrop for our constructions is provided by the
theory of {\em determinant groupoids} over $k$. The notion of
determinant groupoid is a natural outgrowth of the theory of the Sato
Grassmannian and hence is  quite familiar in mathematical physics. But
inconveniently, in its native habitat, this notion is packaged with a
lot of extra structure unneeded for studying reciprocity laws. We sketch
here a ``minimalist'' version of the theory just adequate for the
purposes we have in mind. We have taken pains to make the
theory concrete and easy 
to apply, and also
suitable for study by beginning graduate students in number theory
and algebraic geometry.  The theory very likely has
applications beyond those discussed in this
paper.   We hope
that our approach can be
generalized to yield an ``integrated version'' of Beilinson's
multi-dimensional generalization
\cite{Beilinson} of Tate's theory.

\section{Determinant groupoids over an artinian local 
ring}
We fix an artinian local ring $k$ throughout the paper.
We denote the
maximal ideal of $k$ by $m$ and the multiplicative group of $k$ by
$k^\times$. Given a free
$k$-module
$V$ of finite rank,  we denote the rank of $V$ over $k$
by $\rank V$ and the maximal exterior power of $V$ over $k$ by 
$\det V$.

\subsection{Background on $k$-modules}
\begin{Lemma}\label{Lemma:FlatnessCriterion}
A $k$-module $V$  is flat if and only if for all integers $r>0$,
row vectors $x\in \Mat_{1\times r}(k)$ and 
column vectors
$v\in
\Mat_{r\times 1}(V)$
 such that $xv=0$,
there exists an integer $s>0$,
a matrix $y\in \Mat_{r\times s}(k)$,
and  a column vector $w\in
\Mat_{s\times 1}(V)$ such that
$v=yw$ and $xy=0$.
\end{Lemma}
\proof This is a standard flatness criterion holding over any
commutative ring with unit. See \cite[(3.A), Thm.\ 1, part 6, pp.\
17-18]{Matsumura}.
\qed

\begin{Lemma}\label{Lemma:Generation}
A family $\{v_i\}_{i\in I}$ of elements of a $k$-module $V$ 
generates $V$ over
$k$ if and only if the family $\{v_i\bmod{mV}\}_{i\in I}$
generates $V/mV$ over $k/m$. (This is a version of Nakayama's Lemma.
Note that $V$ need not be finitely generated over $k$.)
\end{Lemma}
\proof ($\Rightarrow$) Trivial. ($\Leftarrow$) Let $V'$
be the $k$-span of $\{v_i\}_{i\in I}$. We have
$$V=V'+mV=V'+m^2V=...=V'$$ because the ideal $m$ is nilpotent. 
\qed

\begin{Lemma}\label{Lemma:Independence}
A family $\{v_i\}_{i\in I}$ of elements of a flat $k$-module $V$
is $k$-linearly independent if and only if
the family $\{v_i\bmod{mV}\}_{i\in I}$
of elements of $V/mV$ is $(k/m)$-linearly
independent.
\end{Lemma}
\proof We may assume that $I=\{1,\dots,r\}$. For
convenience assemble the vectors $v_i$ into a column vector $v\in
\Mat_{r\times 1}(k)$. ($\Rightarrow$) Suppose that there exists $x\in
\Mat_{1\times r}(k)$ such that $x\not\equiv 0\bmod{m}$ but $xv\equiv
0\bmod{mV}$. Let $T$ be a minimal ideal of $k$ and select $0\neq
t\in T$. Then $tm=0$, hence $txv=0$ but $tx\neq 0$, a contradiction.
(Flatness of $V$ was not needed to prove this implication.)
($\Leftarrow$) Suppose there exists $x\in \Mat_{1\times r}(k)$
such that $xv=0$. By Lemma~\ref{Lemma:FlatnessCriterion} there exists
an integer
$s>0$, a matrix
$y\in
\Mat_{r\times s}(k)$ and a  column vector $w\in\Mat_{s\times
1}(k)$ such that $v=yw$ and $xy=0$. By hypothesis the matrix
$y\bmod{m}\in \Mat_{r\times s}(k/m)$ must be of maximal rank.
Therefore some maximal square submatrix of $y$ is invertible
and hence $x=0$. (This argument comes from
the proof of
\cite[(3.G), Proposition, p.\ 21]{Matsumura}.) \qed

\begin{Proposition}
(i) A $k$-module is free if and only if
projective if and only if flat.  
(ii) If two of the $k$-modules in a short exact sequence
of such are free, then so is the third. 
(iii) The $k$-linear dual of a free $k$-module is again free.
(We frequently apply this proposition below but rarely cite it
explicitly.)
\end{Proposition}
\proof (i) It is necessary only to prove that flatness implies
freeness, and for this purpose Lemmas~\ref{Lemma:Generation} and
\ref{Lemma:Independence} clearly suffice. (ii)  Let 
$$0\rightarrow U\rightarrow V\rightarrow W\rightarrow 0$$
be a short exact sequence of $k$-modules. If $U$ and $W$
are free, then $V$ is clearly free. 
 If
$V$ and
$W$ are free then $U$ is a direct summand of $V$, hence $U$ is
projective and hence $U$ is free. If $U$ and $V$ are free,
then any $k$-basis $\{u_i\}$ of $U$ remains $(k/m)$-linearly
independent in $V/mV$ by Lemma~\ref{Lemma:Independence}, hence there
exists a family of elements
$\{v_j\}$ of $V$ such that $\{u_i\bmod{mV}\}\coprod \{v_j\bmod{mV}\}$
is a $(k/m)$-basis of $V/mV$, hence $\{u_i\}\coprod\{v_j\}$
is a $k$-basis of $V$ by
Lemmas~\ref{Lemma:Generation} and \ref{Lemma:Independence}, and
hence
$\{v_j\}$ projects to a $k$-basis for $W$. (iii) This can be proved
by straightforwardly applying Lemma~\ref{Lemma:FlatnessCriterion} to
verify flatness. We omit the details.
\qed

\begin{Lemma}\label{Lemma:TransitivityChecker}
Let $V$ be a free $k$-module and let $M$ be a finitely generated
$k$-module. For every
$k$-linear map
$f:V\rightarrow M$ there exists a $k$-submodule $V'\subset V$
such that the quotient $k$-module $V/V'$ is free of finite rank and
$V'\subset \ker f$. 
\end{Lemma}
\proof By induction on the length of
$M$ as a
$k$-module, we may assume that $M=k/m$ and $f\neq 0$.
Since $V$ is free, $f$ is the reduction modulo $m$
of a surjective $k$-linear functional $\tilde{f}:V\rightarrow
k$. Put $V':=\ker\tilde{f}$. Then $V'$ has the desired properties.
\qed

\begin{Proposition}\label{Proposition:TransitivityChecker}
Let $V$ be a free $k$-module. Let $A,B\subset V$ be free $k$-submodules.
The quotient $(A+B)/(A\cap B)$ is finitely generated if and only if
there exists a free
$k$-submodule $P\subset A\cap B$ such that
the quotient $k$-modules $A/P$ and $B/P$ are free of
finite rank. (We frequently apply
this proposition below but rarely cite it explicitly.)
\end{Proposition}
\proof ($\Rightarrow$) By
Lemma~\ref{Lemma:TransitivityChecker} applied to the natural map
$B\rightarrow B/(A\cap B)$ there exists a $k$-submodule $P\subset A\cap
B$ such that the quotient $k$-module $B/P$ is free of finite rank.
It is then easily verified that $P$ has all the desired properties.
($\Leftarrow$) Trivial. 
\qed

\subsection{Commensurability and related notions}
\subsubsection{Commensurability}
Let $V$ be a free $k$-module. Given free $k$-submodules $A,B\subset V$
we write $A\sim B$ and say that $A$ and $B$ are {\em commensurable}
if the quotient $\frac{A+B}{A\cap B}$ is finitely generated over
$k$.
It is easily verified that commensurability
is an equivalence
relation.

\subsubsection{Relative rank}
 
Let $V$ be a free $k$-module. Let
$A,B\subset V$ be free
$k$-submodules such that $A\sim B$. Put
$$\Rank^V(A,B):=\rank A/P-\rank B/P$$
where $P\subset A\cap B$ is any free $k$-submodule such that
$A/P$ and $B/P$ are free of finite rank, thereby defining the
{\em relative rank} of
$A$ and $B$. It is easily verified that
$\Rank^V(A,B)$ is independent of the choice of $P$
and hence well-defined. It follows that 
$$\Rank^V(A,C)=\Rank^V(A,B)+\Rank^V(B,C)$$
for all free $k$-submodules $A,B,C\subset V$
such that $A\sim B\sim C$.

\subsubsection{The restricted general linear group}
Given a free $k$-module $V$ and a free $k$-submodule $V_+\subset
V$, let
$G^V_{V_+}$ denote the set of
$k$-linear automorphisms $\sigma$ of
$V$ such that $V_+\sim \sigma V_+$. 
It is
easily verified that
$G^V_{V_+}$ is a subgroup of the group of $k$-linear automorphisms of
$V$ depending only on the commensurability class of
$V_+$. We call $G^V_{V_+}$ the {\em restricted general linear
group} associated to $V$ and $V_+$.

\subsubsection{The index} Let $V$ be a free $k$-module
and let 
$V_+\subset V$ be a free
$k$-submodule. Given
$\sigma\in G^V_{V_+}$, put
$$\ind^V_{V_+} \sigma:= \Rank(V_+,\sigma V_+),$$
thereby defining the {\em index} of
$\sigma$.  It is easily verified that $\ind^V_{V_+}\sigma$ depends only
on the commensurability class of $V_+$. It follows that the function
$$\ind^V_{V_+}:G_{V_+}^V\rightarrow \ZZ$$
is a group homomorphism.
\begin{Lemma}\label{Lemma:UpsideDown}
Let $V$ be a free $k$-module equipped with a direct sum decomposition
$V=V_+\oplus V_-$. For all 
$\sigma\in G_{V_-}^V\cap G_{V_+}^V$ we have 
$\ind_{V_+}^V\sigma+\ind_{V_-}^V\sigma=0$.
\end{Lemma}
\proof We have
$$\begin{array}{cl}
&\ind_{V_+}^V\sigma+\ind_{V_-}^V\sigma\\
=&\rank V_+/P_+ -\rank \sigma V_+/P_+ +\rank V_-/P_- -\rank \sigma
V_-/P_-\\ =&\rank (V_++V_-)/(P_++P_-)-\rank(\sigma V_++\sigma
V_-)/(P_++P_-)\\ =&\rank V/(P_++P_-)-\rank V/(P_++P_-)=0
\end{array}
$$
for any free $k$-submodules $P_{\pm}\subset V_{\pm}\cap \sigma V_{\pm}$
such that $V_{\pm}/P_{\pm}$ and $\sigma V_{\pm}/P_{\pm}$
are free of finite rank. \qed
\pagebreak

\subsection{The construction $\Det^V(A,B;P)$.}
 Fix a free $k$-module $V$.
\subsubsection{Definition}
Given commensurable free $k$-submodules
$A,B\subset V$
and a free $k$-submodule $P\subset A\cap B$ such that the quotient
$k$-modules
$A/P$ and $B/P$ are free of finite rank, put
$$ \Det^V(A,B;P):=
\left\{\mbox{$k$-linear isomorphisms $\det(A/P)\iso
\det(B/P)$}\right\}.$$
Note that we have at our disposal an operation of
{\em scalar multiplication}
$$((x,\alpha)\mapsto x\cdot\alpha):
k^\times \times \Det^V(A,B;P)\rightarrow \Det^V(A,B;P)$$
with respect to which $\Det^V(A,B;P)$ becomes
a $k^\times$-torsor.

\subsubsection{The composition law}
Given free $k$-submodules $A,B,C\subset V$ belonging to the same
commensurability class and a free $k$-submodule $P\subset A\cap B\cap
C$  such that the $k$-modules $A/P$, $B/P$ and $C/P$ are free of finite
rank, we have a {\em composition law}
$$
((\alpha,\beta)\mapsto \beta\circ \alpha):
\Det^V(A,B;P)\times \Det^V(B,C;P)\rightarrow
\Det^V(A,C;P)$$
at our disposal.
The composition law is compatible with scalar multiplication
in the sense that 
$$x\cdot(\beta\circ \alpha)=(x\cdot\beta)\circ
\alpha=\beta\circ(x\cdot\alpha).$$ Clearly the composition law is
associative.

\subsubsection{The cancellation rule}
Given commensurable free
$k$-submodules $A,B\subset V$
and free $k$-submodules $Q\subset P\subset A\cap B$
such that the quotient $k$-modules $A/Q$ and $B/Q$ are free of finite
rank,  we
have at our disposal a canonical isomorphism 
$$\left(\left(\wedge p_\ell\right)\wedge
\left(\wedge
\tilde{a}_i\right)\mapsto \left(\wedge p_\ell\right)\wedge
\left(\wedge \tilde{b}_j\right)\right)
\mapsto
\left(
\left(\wedge
a_i\right)\mapsto 
\left(\wedge b_j\right)\right):$$
$$\Det^V(A,B;Q)\iso
\Det^V(A,B;P)$$
where $\{a_i\}$, $\{b_j\}$ and $\{p_\ell\}$
are any ordered $k$-bases for $A/P$, $B/P$ and $P/Q$, respectively,
and $\{\tilde{a}_i\}$ and $\{\tilde{b}_j\}$ are liftings
to $A/Q$ and $B/Q$, respectively, of $\{a_i\}$ and $\{b_j\}$,
respectively. We refer to this isomorphism as the {\em cancellation
rule}.  It is easily verified that the
cancellation rule commutes with the composition law and with scalar
multiplication.

\subsubsection{Inverse system structure}
\label{subsubsection:InverseSystem}
Let $V$ be a free $k$-module and let $A,B\subset V$ be commensurable
free $k$-submodules. Consider the family
$\PP$ of free
$k$-submodules
$P\subset A\cap B$ such that $A/P$ and $B/P$ are free of finite rank.
Partially order the  family $\PP$ by reverse inclusion.
Then 
$\PP$ is a directed set and it is easily
verified that the cancellation rule gives the family of sets
$\Det^V(A,B;P)$ indexed by 
$P\in \PP$  the structure of inverse system. 

\pagebreak

\subsection{The connected groupoid $\Det^V_{V_+}$}
Fix a free $k$-module $V$ and a free $k$-submodule $V_+\subset V$.
\subsubsection{Definition}
For all free $k$-submodules $A,B\subset V$
belonging to the commensurability class of $V_+$
put
$$\Det^V_{V_+}(A,B):=
\lim_{\leftarrow}\Det^V(A,B;P),
$$
the limit extended over the inverse system defined in
\S\ref{subsubsection:InverseSystem}. Since scalar multiplication
commutes with the cancellation rule, we have at our disposal an
operation
 of scalar multiplication
$$\left(\left(x,\alpha\right)\mapsto x\cdot\alpha\right):
k^\times\times\Det^V_{V_+}(A,B)\rightarrow
\Det^V_{V_+}(A,B)$$
endowing the set $\Det^V_{V_+}(A,B)$ with the
structure of $k^\times$-torsor.

\subsubsection{The composition law}
Since the cancellation rule and the composition law 
commute, we obtain in the limit a composition law
$$((\alpha,\beta)\mapsto \beta\circ \alpha):\Det^V_{V_+}(A,B)
\times \Det^V_{V_+}(B,C)\rightarrow
\Det^V_{V_+}(A,C)$$
for all free $k$-submodules $A,B,C\subset V$ belonging to the same
commensurability class. The composition law is associative and moreover
is compatible with scalar multiplication in the evident sense.

\subsubsection{Definition of $\Det^V_{V_+}$}
The rule sending each pair $A,B\subset V$ of
free $k$-submodules commensurable to $V_+$ to the set
$\Det_{V_+}^V(A,B)$ makes the commensurability class of $V_+$ into a
category. We denote this category by
$\Det^V_{V_+}$. It is easily verified that every morphism
in $\Det^V_{V_+}$ is an isomorphism and that all the objects of
$\Det^V_{V_+}$ are isomorphic.
Thus $\Det^V_{V_+}$ is a connected groupoid. Note that
$\Det^V_{V_+}$ depends only on the commensurability class of $V_+$.

\subsubsection{Abstract nonsense}
Fix $\sigma\in G_{V_+}^V$ and 
let
$\sigma_*$ be the functor from $\Det^V_{V_+}$ to itself
induced in evident fashion by $\sigma$.
A {\em natural transformation} $\Phi$ from the identity
functor of $\Det^V_{V_+}$ to $\sigma_*$ is a rule
associating to each object $A$ of $\Det^V_{V_+}$ 
a morphism $\Phi(A)\in
\Det_{V_+}^V(A,\gamma A)$ such that for any two
objects $A$ and $B$ of $\Det^V_{V_+}$ the diagram
$$\begin{array}{rcccl}
&A&\xrightarrow{\Phi(A)}&\gamma A\\
\alpha&\downarrow&&\downarrow&\gamma_*\alpha\\
&B&\xrightarrow{\Phi(B)}&\gamma B
\end{array}
$$
commutes for all $\alpha\in \Det_{V_+}^V(A,B)$.
It is easily verified that the natural transformations
from the identity functor of $\Det^V_{V_+}$ to $\sigma_*$
are in bijective correspondence with $\Det^V_{V_+}(V_+,\gamma V_+)$
via the map $\Phi\mapsto \Phi(V_+)$. 

\subsubsection{The central extension
$\tilde{G}^V_{V_+}$} We define $\tilde{G}^V_{V_+}$ to be the set
consisting of pairs $(\sigma,\Phi)$ where
$\sigma\in G^V_{V_+}$ and
$\Phi$ is a natural transformation from the identity functor of
$\Det^V_{V_+}$ to $\sigma_*$. We compose elements of $\tilde{G}^V_{V_+}$
by the rule
$$(\sigma_1,\Phi_1)(\sigma_2,\Phi_2):=
(\sigma_1\sigma_2, A\mapsto
(\sigma_{1*}(\Phi_2(A)))\circ\Phi_1(A)).$$
It is easily verified that this composition law is a group law.
The group
$\tilde{G}^V_{V_+}$ thus defined depends only on the
commensurability class of
$V_+$.  The group $\tilde{G}^V_{V_+}$ fits into a
canonical exact sequence
$$1\rightarrow k^\times \xrightarrow{x\mapsto(1,A\mapsto x\cdot 1_A)}
\tilde{G}^V_{V_+}\xrightarrow{(\sigma,\Phi)\mapsto \sigma}
G^V_{V_+}\rightarrow 1,$$
where $1_A\in \Det^V_{V_+}(A,A)$ denotes the identity map.
The exact sequence identifies $k^\times$ with a subgroup
of the center of $\tilde{G}^V_{V_+}$. 

\subsubsection{Remark}
When $k$ is a field, the central extension of the restricted
general linear group $G^V_{V_+}$  constructed here has cohomology
class in $H^2(G^V_{V_+},k^\times)$ equal to the
cohomology class associated to the central extension studied in
the paper \cite{AdCK} and equal to the opposite of the
cohomology class associated to the central extension studied in
\cite{PablosRomo}.
 So what we study here should be regarded
as the deformation theory of the central extensions
of \cite{AdCK} and \cite{PablosRomo}.

\section{Study of the symbol $\{\sigma,\tau\}_{V_+}^V$}
\subsection{Definition and basic properties of the symbol}
Fix a free $k$-module $V$ and free $k$-submodule $V_+\subset V$.
\subsubsection{Definition}\label{subsubsection:SymbolDefinition}
Given
commuting elements
$\sigma,\tau\in G^V_{V_+}$, put
$$\{\sigma,\tau\}_{V_+}^V:=\tilde{\sigma}\tilde{\tau}\tilde{\sigma}^{-1}\tilde{\tau}^{-1}
\in \ker\left(\tilde{G}^V_{V_+}\rightarrow
G^V_{V_+}\right)=k^\times$$
where $\tilde{\sigma},\tilde{\tau}\in
\tilde{G}^V_{V_+}$ are any liftings of $\sigma$ and $\tau$,
respectively. It is easily verified that $\{\sigma,\tau\}_{V_+}^V$ is
independent of the choice of liftings and hence well-defined. 
By the definitions we have
$$(\sigma_*\beta)\circ \alpha=\{\sigma,\tau\}_{V_+}^V
\cdot((\tau_*\alpha)\circ \beta)$$
for all $\alpha\in \Det^V_{V_+}(V_+,\sigma V_+)$
and $\beta\in \Det^V_{V_+}(V_+,\tau V_+)$. The latter
formula is the one we rely upon in practice to calculate
$\{\sigma,\tau\}_{V_+}^V$.

\begin{Lemma}
Let $G$ be a group. Write $[x,y]:=xyx^{-1}y^{-1}$
for all $x,y\in G$. Now let
$\alpha,\beta,\gamma\in G$ be given.
Assume that $[\beta,\gamma]$ is central in $G$. Then we have
$[\alpha,\gamma][\beta,\gamma]=[\alpha\beta,\gamma]$.
\end{Lemma}
\proof 
$$\begin{array}{rcl}
[\alpha,\gamma][\beta,\gamma][\alpha\beta,\gamma]^{-1}&
=&\alpha \gamma \alpha^{-1}\gamma^{-1}[\beta,\gamma]\gamma
\alpha\beta\gamma^{-1}\beta^{-1}\alpha^{-1}\\
&=&\alpha \gamma
[\beta,\gamma]\beta\gamma^{-1}\beta^{-1}\alpha^{-1}\\
&=&
\alpha [\beta,\gamma][\beta,\gamma]^{-1}\alpha^{-1}
=1.
\end{array}
$$
\qed

\subsubsection{Basic properties}
\label{subsubsection:BasicProperties}
Fix elements $\sigma,\sigma',\tau,\tau'\in G^V_{V_+}$
such that the $\sigma$'s commute with the $\tau$'s.
(But we need assume neither that
$\sigma\sigma'=\sigma'\sigma$ nor that $\tau\tau'=\tau'\tau$.) 
The
following relations hold:\\

\begin{itemize}
\item $\{\sigma,\sigma\}_{V_+}^V=1$.\\
\item
$\{\sigma,\tau\}_{V_+}^V=\left(\{\tau,\sigma\}_{V_+}^V\right)^{-1}$.\\
\item
$\{\sigma\sigma',\tau\}_{V_+}^V
=\{\sigma,\tau\}_{V_+}^V\{\sigma',\tau\}_{V_+}^V$.\\
\item
$\{\sigma,\tau\tau'\}_{V_+}^V
=\{\sigma,\tau\}_{V_+}^V\{\sigma,\tau'\}_{V_+}^V$.\\
\item $\sigma V_+=V_+=\tau V_+\Rightarrow
\{\sigma,\tau\}_{V_+}^V=1$.\\
\item If $V_+=\{0\}$
or
$V_+=V$, then $\{\sigma,\tau\}_{V_+}^V=1$.\\
\item 
$\{\sigma,\tau\}_{V_+}^V$ depends only on the commensurability class of
$V_+$.\\
\end{itemize}
For the most part these facts are straightforwardly deduced from the
definitions. Only the proof of bimultiplicativity offers any
difficulty and  
the essential point of that proof is contained in the preceding
lemma.

\subsection{The four square
identity}\label{subsection:FourSquareIdentity} Fix a free $k$-module
$V$, a free $k$-submodule $V_+\subset V$ and commuting elements
$\sigma,\tau\in G^V_{V_+}$.  We work out an explicit formula for the
symbol 
$\{\sigma,\tau\}_{V_+}^V$ in terms of determinants.

\subsubsection{Choices}
Let 
$$P\subset V_+\cap \sigma V_+,\;\;\;Q\subset V_+\cap \tau V_+$$ be
free
$k$-submodules such that the quotient $k$-modules
$$V_+/P,\;\;\;\sigma V_+/P,\;\;\;V_+/Q,\;\;\;\tau V_+/Q$$
are free of finite rank.
Let 
$$R\subset P\cap Q\cap \tau P\cap \sigma Q$$
be a free $k$-submodule such that the quotient
$k$-modules 
$$P/R,\;\;\;Q/R,\;\;\;\tau P/R,\;\;\;\sigma Q/R$$
are free of finite rank. Fix finite sequences
$$\{a_i\},\dots, \{h_i\}$$
in $V$ (not necessarily all of the same length) such that the
corresponding finite sequences
$$
\left.
\begin{array}{l}
\{a_i\bmod{P}\}, \\
\{b_i\bmod{P}\}, \\
\{c_i\bmod{Q}\}, \\
\{d_i\bmod{Q}\},\\
\{e_i\bmod{R}\},\\
\{f_i\bmod{R}\},\\
\{g_i\bmod{R}\},\\
\{h_i\bmod{R}\},
\end{array}\right\}\;\;\mbox{are $k$-bases of}\;\;
\left\{\begin{array}{l}
V_+/P,\\
\sigma
V_+/P,\\
V_+/Q,\\
\tau_+V/Q,\\
P/R,\\
Q/R,\\
\tau P/R,\\
\sigma Q/R,\\
\end{array}\right.
$$
respectively.
Fix morphisms
$$\alpha\in \Det^V_{V+}(V_+,\sigma V_+),\;\;\;
\beta\in \Det^V_{V_+}(V_+,\tau V_+).$$
\subsubsection{Construction of representations of the morphisms
$\alpha$,
$\beta$}
Since our purpose is to calculate $\{\sigma,\tau\}_{V_+}^V$ by means of
the last formula of \S\ref{subsubsection:SymbolDefinition}, we may simply 
assume that
$\alpha$ is represented by 
$$\left(\wedge (a_i\bmod{P})\mapsto \wedge (b_i\bmod{P})\right)
\in \Hom_k\left(\det(V_+/P),\det(\sigma V_+/P)\right)$$
and that $\beta$ is represented by
$$\left(\wedge (c_i\bmod{Q})\mapsto \wedge (d_i\bmod{Q})\right)
\in \Hom_k\left(\det(V_+/Q),\det(\tau V_+/Q)\right).$$
By the cancellation rule $\alpha$ is also represented by
$$\left(\left(\wedge \bar{e}_i\right)
\wedge \left(\wedge \bar{a}_i\right)\mapsto 
\left(\wedge \bar{e}_i\right)
\wedge \left(\wedge \bar{b}_i\right)\right)
\in \Hom_k\left(\det(V_+/R),\det(\sigma V_+/R)\right)
$$
and $\beta$ is also represented by
$$\left(\left(\wedge \bar{f}_i\right)
\wedge \left(\wedge \bar{c}_i\right)\mapsto 
\left(\wedge \bar{f}_i\right)
\wedge \left(\wedge \bar{d}_i\right)\right)
\in \Hom_k\left(\det(V_+/R),\det(\tau V_+/R)\right),$$
where here and below $v\mapsto \bar{v}$
denotes reduction modulo $R$. 

\subsubsection{Construction of representations of the morphisms
$\tau_*\alpha$ and $\sigma_*\beta$}
By definition $\tau_*\alpha$ is represented by 
$$\begin{array}{l}
\left(\wedge (\tau a_i\bmod{\tau P})\mapsto
\wedge (\tau b_i\bmod{\tau P})\right)\\\\
\in \Hom_k\left(\det(\tau V_+/\tau P),\det(\tau\sigma
V_+/\tau P)\right)
\end{array}
$$ and $\sigma_*\beta$ is represented by
$$\begin{array}{l}
\displaystyle\left(\wedge (\sigma c_i\bmod{\sigma Q})\mapsto
\wedge (\sigma d_i\bmod{\sigma Q})\right)\\\\
\in \Hom_k\left(\det(\sigma V_+/\sigma Q),\det(\sigma\tau
V_+/\sigma Q)\right).
\end{array}$$
By the cancellation rule, $\tau_*\alpha$ is also represented by
$$\begin{array}{l}
\displaystyle\left(\left(\wedge \bar{g}_i\right)\wedge
\left(\wedge \overline{\tau a}_i\right)
\mapsto \left(\wedge \bar{g}_i\right)
\wedge \left(\wedge \overline{\tau b}_i\right)\right)\\\\
\in \Hom_k\left(\det(\tau V_+/R),\det(\tau\sigma V_+/R)\right)
\end{array}$$
and $\sigma_*\beta$ is also represented by
$$\begin{array}{l}
\displaystyle\left(\left(\wedge \bar{h}_i\right)\wedge
\left(\wedge \overline{\sigma c}_i\right)
\mapsto \left(\wedge \bar{h}_i\right)
\wedge \left(\wedge \overline{\sigma d}_i\right)\right)\\\\
\in \Hom_k\left(\det(\sigma V_+/R),\det(\sigma \tau V_+/R)\right).
\end{array}$$

\subsubsection{Conclusion of the calculation}
We have now represented all the morphisms $\alpha$, $\beta$,
$\tau_*\alpha$ and $\sigma_*\beta$ in such a way that we can
obtain representations for the compositions $(\sigma_*\beta)\circ
\alpha$ and $(\tau_*\alpha)\circ \beta$ in the same
rank one free $k$-module, namely
$$\Hom_k\left(\det(V_+/R),\det(\sigma \tau V_+/R)\right)
=\Hom_k\left(\det(V_+/R),\det(\tau\sigma V_+/R)\right).$$
It is a straightforward matter to calculate the ratio.
We find that
$$\begin{array}{cl}
&\{\sigma,\tau\}_{V_+}^V\\\\
=&\displaystyle\frac{(\wedge \bar{h}_i)\wedge (\wedge
\overline{\sigma d}_i)}
{(\wedge \bar{g}_i)\wedge (\wedge \overline{\tau b}_i)}
\frac{(\wedge \bar{e}_i)\wedge (\wedge \bar{b}_i)}
{(\wedge \bar{h}_i)\wedge (\wedge \overline{\sigma c}_i)}
\frac{(\wedge \bar{f}_i)\wedge (\wedge \bar{c}_i)}
{(\wedge \bar{e}_i)\wedge (\wedge \bar{a}_i)}
\frac{(\wedge \bar{g}_i)\wedge (\wedge \overline{\tau a}_i)}
{(\wedge \bar{f}_i)\wedge (\wedge \bar{d}_i)}.
\end{array}
$$
The right side of the formula makes sense because in each of the
fractions, both
numerator and denominator generate the same rank one $k$-submodule of
the exterior algebra over $k$ of $V/R$.  We call this result the {\em
four square identity} because in an
obvious way the diagram\\
$$\begin{array}{cccccc}
V_+&\{a_i\}&P&\{b_i\}&\sigma V_+\\\\
\{c_i\}&&\{e_i\}&&\{\sigma c_i\}\\\\
Q&\{f_i\}&R&\{h_i\}&\sigma Q\\\\
\{d_i\}&&\{g_i\}&&\{\sigma d_i\}\\\\
\tau V_+&\{\tau a_i\}&\tau P&\{\tau b_i\}&\sigma \tau V_+
&=\tau\sigma V_+
\end{array}$$\\
serves as a mnemonic. We call the diagram above a {\em template} for
the calculation of
$\{\sigma,\tau\}^V_{V_+}$ and we say that the right side of the four
square identity is the {\em value} of the template.

\subsection{General rules of calculation}
Fix a free $k$-module $V$ and a free $k$-submodule $V_+\subset V$.

\begin{Proposition}\label{Proposition:ExpansionAndContraction}
Fix commuting elements $\sigma,\tau\in G_{V_+}^V$.
(i) Suppose there exists a free $k$-module $W$ containing
$V$ as a $k$-submodule. Suppose further that $\sigma,\tau$
admit commuting extensions $\tilde{\sigma},\tilde{\tau}\in G^W_{V_+}$,
respectively. Then we have
$\{\sigma,\tau\}_{V_+}^V=\{\tilde{\sigma},\tilde{\tau}\}_{V_+}^W$.
(ii) Suppose there exists a free $k$-submodule $U\subset V$ such that
$\sigma U=U=\tau U$ and $U\cap V_+=0$. Put $\bar{V}:=V/U$ and
$\bar{V}_+:=(V_++U)/U$. Let
$\bar{\sigma}$ and
$\bar{\tau}$ be the $k$-linear automorphisms of $\bar{V}$ induced by
$\sigma$ and $\tau$, respectively, and assume that
$\bar{\sigma},\bar{\tau}\in G^{\bar{V}}_{\bar{V}_+}$. Then we have
$\{\sigma,\tau\}_{V_+}^V=
\{\bar{\sigma},\bar{\tau}\}^{\bar{V}}_{\bar{V}_+}.$
\end{Proposition}
\proof We return to the setting of
\S\ref{subsection:FourSquareIdentity}. Let
$T$ be the template above for calculating
$\{\sigma,\tau\}_{V_+}^V$. The very same template $T$
serves also to calculate $\{\tilde{\sigma},\tilde{\tau}\}_{V_+}^W$.
Therefore (i) holds. Let $\bar{T}$ be the projection of the template $T$
into $\bar{V}$. Then $\bar{T}$ is a template for the
calculation of
$\{\bar{\sigma},\bar{\tau}\}^{\bar{V}}_{\bar{V}_+}$ and
moreover it is easily verified that the values of the templates $T$ and
$\bar{T}$ are equal. Therefore (ii) holds. \qed

\begin{Proposition}
\label{Proposition:Flab}
Suppose $V$ is equipped with a direct sum
decomposition
$V=V_0\oplus V_1$. Put $V_{i+}:=V_i\cap V_+$ for $i=0,1$ and assume that
$V_+=V_{0+}\oplus V_{1+}$. Let commuting elements $\sigma_0,\sigma_1\in G^V_{V_+}$ be
given such that
$$\sigma_i\vert_{V_0}\in G_{V_{0+}}^{V_0},\;\;\;
\sigma_i\vert_{V_1}=1$$
for $i=0,1$. Then we have
$$\left\{\sigma_0\vert_{V_0},\sigma_1\vert_{V_0}\right\}_{V_{0+}}^{V_0}
=\{\sigma_0,\sigma_1\}_{V_+}^V.$$
\end{Proposition}
\proof The proof is a straightforward application of the four square
identity similar to that made in the proof of
Proposition~\ref{Proposition:ExpansionAndContraction} and therefore
omitted.
\qed

\begin{Proposition}\label{Proposition:SignRule}
Again suppose $V$ is equipped with a direct sum decomposition
$V=V_0\oplus V_1$, put $V_{i+}:=V_i\cap V_+$ for $i=0,1$ and assume that
$V_+=V_{0+}\oplus V_{1+}$. Let $\sigma_0,\sigma_1\in G^V_{V_+}$ be
given such that
$$\sigma_i\vert_{V_i}\in
G^{V_i}_{V_{i+}},\;\;\;\sigma_i\vert_{V_{1-i}}=1$$ for $i=0,1$.
(Necessarily $\sigma_0$ and $\sigma_1$ commute.) Then we
have
$$\{\sigma_0,\sigma_1\}_{V_+}^V=
(-1)^{\nu_0\nu_1}
$$
where
$$\nu_i:=\ind_{V_{i+}}^{V_i}\sigma_i\vert_{V_i}=
\ind_{V_+}^V\sigma_i$$
for $i=0,1$.
\proof
For $i=0,1$, choose a free $k$-submodule $P_i\subset V_{i+}\cap
\sigma_i V_{i+}$ such that 
the quotient $k$-modules $V_{i+}/P$ and $\sigma_iV_{i+}/P_i$ are free
of finite rank and  also choose  finite sequences $\{e_{ij}\}$ and
$\{f_{ij}\}$  in $V_i$ (not necessarily of the same
length) reducing modulo
$P_i$ to $k$-bases for
$V_{i+}/P_i$ and $\sigma_iV_{i+}/P_i$, respectively.
Then the diagram\\
$$
\begin{array}{ccccc}
V_{0+}\oplus V_{1+}&\{e_{0j}\}&P_0\oplus V_{1+}&
\{f_{0j}\}&\sigma_0V_{0+}\oplus V_{1+}\\\\
\{e_{1j}\}&&\{e_{1j}\}&&\{e_{1j}\}\\\\
V_{0+}\oplus P_1&\{e_{0j}\}&P_0\oplus
P_1&\{f_{0j}\}&\sigma_0V_{0+}\oplus P_1\\\\
\{f_{1j}\}&&\{f_{1j}\}&&\{f_{1j}\}\\\\
V_{0+}\oplus \sigma_1V_{1+}&\{e_{0j}\}&P_0\oplus \sigma_1V_{1+}
&\{f_{0j}\}&\sigma_0V_{0+}\oplus \sigma_1 V_{1+}
\end{array}
$$\\
is a template for the calculation of
$\{\sigma_0,\sigma_1\}_{V_+}^V$. The desired result now follows 
by the four square identity and the definitions.
\qed
\end{Proposition}
\begin{Proposition}\label{Proposition:Switcheroo}
Let $V_-\subset V$ be a free $k$-submodule such that $V=V_+\oplus
V_-$.  Let commuting elements
$\sigma,\tau\in G^V_{V_+}\cap G^V_{V_-}$ be given. Then we have
$$\{\sigma,\tau\}^V_{V_+}\{\sigma,\tau\}^V_{V_-}=1.$$
\end{Proposition}
\proof
We define
$$\sigma_0,\sigma_1,\tau_0,\tau_1\in G^{V\oplus V}_{V_+\oplus V_-}$$
by the block decompositions
$$\sigma_0=\left[\begin{array}{cc}
\sigma&0\\
0&1\end{array}\right],\;\;\;
\sigma_1=\left[\begin{array}{cc}
1&0\\
0&\sigma\end{array}\right],\;\;\;
\tau_0=\left[\begin{array}{cc}
\tau&0\\
0&1\end{array}\right],\;\;\;
\tau_1=\left[\begin{array}{cc}
1&0\\
0&\tau\end{array}\right].
$$
Then we have\\
$$\begin{array}{cl}
&\{\sigma_0\sigma_1,\tau_0\tau_1\}_{V_+\oplus
V_-}^{V\oplus V}\\\\
 =&\displaystyle\prod_{i=0}^1\prod_{j=0}^1
\{\sigma_i,\tau_j\}_{V_+\oplus V_-}^{V\oplus V}=
\{\sigma,\tau\}_{V_+}^V\{\sigma,\tau\}_{V_-}^V
(-1)^{\mu_-\nu_++\mu_+\nu_-}
\end{array}$$\\
by Propositions~\ref{Proposition:Flab}
and \ref{Proposition:SignRule}, where
$$\mu_{\pm}=\ind_{V_{\pm}}^V \sigma,\;\;\;
\nu_{\pm}=\ind_{V_{\pm}}^V\tau.$$
Put
$$U:=\ker
\left(\left(v_0\oplus v_1\mapsto v_0+v_1\right):V\oplus V\rightarrow
V\right).$$ Further, we have\\
$$\begin{array}{rcl}
\{\sigma_0\sigma_1,\tau_0\tau_1\}_{V_+\oplus
V_-}^{V\oplus V}
&=&\{\sigma_0\sigma_1\bmod{U},\tau_0\tau_1\bmod{U}\}_{((V_+\oplus
V_-)+U)/U}^{(V\oplus V)/U}\\\\
&=&\{\sigma,\tau\}^V_V=1
\end{array}$$\\ by part (ii) of
Proposition~\ref{Proposition:ExpansionAndContraction}. Finally,
we have
$$(-1)^{\mu_-\nu_++\mu_+\nu_-}=(-1)^{2\mu_-\nu_+}=1$$
by Lemma~\ref{Lemma:UpsideDown}. The result follows.
\qed

\subsection{The commutator interpretation of the Contou-Carr\`{e}re
symbol}

\subsubsection{Preliminary discussion of the ring $k((t))$}
Let $t$ be a variable. Let $k((t))$ be the ring obtained from the
power series ring
$k[[t]]$ by inverting $t$. It is easily verified that $k((t))$ is an
artinian local ring with maximal ideal $m((t))$ and residue field
$(k/m)((t))$. We have an additive direct sum decomposition
$$k((t))=t^{-1}k[t^{-1}]\oplus k[[t]]$$
and a multiplicative direct sum
decomposition
$$k((t))^\times=t^\ZZ\cdot (1+m[t^{-1}])\cdot k^\times \cdot
(1+tk[[t]])$$
at our disposal. The latter decomposition can be refined as follows.
Each $f\in k((t))^\times$ has a unique presentation
$$f=t^{w(f)}\cdot a_0\cdot\prod_{i=1}^\infty
\left(1-a_{-i}t^{-i}\right)
\cdot \prod_{i=1}^\infty \left(1-a_it^i\right)
$$
where 
$$w(f)\in \ZZ,\;\;\;a_0\in k^\times,\;\;\;
\left\{\begin{array}{ll}
a_i=0&\mbox{if $i\ll 0$,}\\a_i\in m&\mbox{if $i<0$,}\\
a_i\in k^\times&\mbox{if $i=0$,}\\
a_i\in k&\mbox{if $i>0$.}
\end{array}\right.$$
We call $w(f)$ the {\em winding number} of $f$
and we call $\{a_i\}_{i=-\infty}^\infty$ the family
of {\em Witt parameters} of $f$.
Now view $f$ as a $k$-linear automorphism of the free
$k$-module
$k((t))$. It is easily verified that
$f\in G_{k[[t]]}^{k((t))}$
and that
$$\ind_{k[[t]]}^{k((t))}f=\ind_{(k/m)[[t]]}^{(k/m)((t))}f=
-w(f\bmod{m})=-w(f).$$

\subsubsection{Definition of the Contou-Carr\`{e}re symbol}
Let $f,g\in k((t))^\times$ be given.
Let
$\{a_i\}$ and $\{b_j\}$ be the systems of Witt parameters
associated to $f$ and $g$, respectively.
Put
$$\langle f,g\rangle:=(-1)^{w(f)w(g)}\cdot
\frac{a_0^{w(g)}}{b_0^{w(f)}}
\cdot
\frac{\prod_{i=1}^\infty \prod_{j=1}^\infty
\left(1-a_i^{j/(i,j)}b_{-j}^{i/(i,j)}\right)^{(i,j)}}
{\prod_{i=1}^\infty \prod_{j=1}^\infty
\left(1-a_{-i}^{j/(i,j)}b_{j}^{i/(i,j)}\right)^{(i,j)}}.
$$
The right side of the definition makes sense because all but
finitely many factors in the infinite products differ from $1$.
This definition is due
to Contou-Carr\`{e}re
\cite{ContouCarrere}. We call the map
$$\langle \cdot,\cdot\rangle:k((t))^\times
\times k((t))^\times \rightarrow k^\times$$
defined by the formula above the {\em Contou-Carr\`{e}re symbol}.
The symbol is clearly anti-symmetric:
$$\langle f,g\rangle=\langle g,f\rangle^{-1}.$$ 
Although it is not immediately evident from the definition,
the symbol is also bimultiplicative:
$$\langle f',g\rangle=\langle f,g\rangle\langle f',g\rangle,\;\;\;
\langle f,gg'\rangle=\langle f,g\rangle\langle f,g'\rangle.$$
The following result establishes bimultiplicativity 
of the Contou-Carr\`{e}re symbol as a byproduct.

\begin{Theorem}\label{Theorem:CommutatorInterpretation}
For all $f,g\in k((t))^\times$ we have
$$\langle f,g\rangle^{-1}=(-1)^{w(f)w(g)}\{f,g\}_{k[[t]]}^{k((t))}$$
where on the right side we view $f$ and $g$ as elements
of the restricted general linear group $G^{k((t))}_{k[[t]]}$.
\end{Theorem}

Before beginning the proof proper we prove a couple
of lemmas.  We say that a polynomial
$f\in k[t]$ is {\em distinguished} if $f$ is monic in $t$ 
and $f\equiv t^{\deg f}\bmod{m}$, in which case
necessarily
$w(f)=\deg f$.

\begin{Lemma}\label{Lemma:CrucialDeterminant} Fix a distinguished
polynomial
$f\in k[t]$ of degree $n$ and $g\in k[[t]]^\times$.  (i) The
quotient of the power series ring $k[[t]]$ by its principal ideal
$(f)$ is free over
$k$  and the monomials $1,t,\dots,t^{n-1}$ form a $k$-basis.
(ii) 
We have 
$$\{f,g\}_{k[[t]]}^{k((t))}=\det(g\vert k[[t]]/(f)).$$
\end{Lemma}
\proof Statement (i) is a special case of the Weierstrass Division
Theorem. From statement (i) it follows that the diagram\\
$$\begin{array}{ccccc}
k[[t]]&\emptyset&k[[t]]&\emptyset&k[[t]]\\\\
\emptyset&&\emptyset&&\emptyset\\\\
k[[t]]&\emptyset&k[[t]]&\emptyset&k[[t]]\\\\
\{f^{-1}t^i\}_{i=0}^{n-1}&&\{f^{-1}t^i\}_{i=0}^{n-1}&&
\{gf^{-1}t^i\}_{i=0}^{n-1}\\\\
f^{-1}\cdot k[[t]]&\emptyset&f^{-1}\cdot k[[t]]&\emptyset&f^{-1}\cdot
k[[t]]
\end{array}
$$\\
is a template for the calculation of $\{g,f^{-1}\}_{k[[t]]}^{k((t))}$.
Statement (ii) now follows from the four square identity.
\qed

\begin{Lemma}\label{Lemma:Distinguished}
Let $f,g\in k[t]$ be distinguished polynomials.
We have 
$$\{f,g\}_{k[[t]]}^{k((t))}=1.$$
\end{Lemma}
 \proof Since
$$f,g\in G^{k((t))}_{k[[t]]}\cap
G^{k((t))}_{t^{-1}k[t^{-1}]},$$ we have
$$\{f,g\}_{k[[t]]}^{k((t))}\{f,g\}_{t^{-1}k[t^{-1}]}^{k((t))}=1$$
by Proposition~\ref{Proposition:Switcheroo}.
Let $p$ and $q$ be the degrees of $f$ and $g$, respectively.
Now consider the template\\
$$
\begin{array}{ccccc}
t^{-1}k[t^{-1}]&\emptyset&t^{-1}k[t^{-1}]&\{t^i\}_{i=0}^{p-1}&t^{p-1}k[t^{-1}]\\\\
\emptyset&&\emptyset&&\emptyset\\\\
t^{-1}k[t^{-1}]&\emptyset&t^{-1}k[t^{-1}]&
\{t^i\}_{i=0}^{p-1}&t^{p-1}k[t^{-1}]\\\\
\{t^i\}_{i=0}^{q-1}&&\{t^i\}_{i=0}^{q-1}&&\{ft^i\}_{i=0}^{q-1}\\\\
t^{q-1}k[t^{-1}]&\emptyset&t^{q-1}k[t^{-1}]&\{gt^i\}_{i=0}^{p-1}&
t^{p+q-1}k[t^{-1}]
\end{array}
$$\\
for the calculation of $\{f,g\}_{t^{-1}k[t^{-1}]}^{k((t))}$.
By the four square identity we conclude that the latter symbol
equals unity, and we are done. \qed

\begin{Lemma}\label{Lemma:Clear}
For all $f,g\in k[[t]]^\times$ we have
$$\{f,g\}_{k[[t]]}^{k((t))}=1.$$
\end{Lemma}
\proof In view of the formal properties 
noted in \S\ref{subsubsection:BasicProperties}, this is clear.
\qed

\subsubsection{Proof of Theorem~\ref{Theorem:CommutatorInterpretation}}
Every element of $k((t))^\times$ factors as a power of $t$ times a
distinguished polynomial times a unit of $k[[t]]$.
So after making the evident reductions based upon
Lemmas~\ref{Lemma:Distinguished} and \ref{Lemma:Clear},
we may assume without loss of generality that $f$ is a distinguished
polynomial and that $g\in k[[t]]^\times$.
Moreover, we may assume without loss of generality that $f$ takes the
special form
$t^p-a$ for some positive integer $p$ and $a\in m$.
By Lemma~\ref{Lemma:CrucialDeterminant} we have
$$\{f,g\}_{k[[t]]}^{k((1/t))}=\det(g\vert k[[t]]/(t^p-a)).$$
This justifies the
further assumption without loss of generality 
 that $g(0)=1$.
Now $t$ operates nilpotently on the quotient $k[[t]]/(t^p-a)$,
hence $t^N\equiv 0\bmod{(t^p-a)}$ for some positive integer $N$,
and hence 
$$\det(1+t^Nh\vert k[[t]]/(t^p-a))=1$$ for
all $h\in k[[t]]$.  This justifies the further assumption without loss
of generality that $g=1-bt^q$ for some positive
integer $q$ and $b\in k$. Finally, we have
$$\det\left(1-bt^q\vert k[[t]]/(t^p-a)\right)
=\left(1-a^{q/(p,q)}b^{p/(p,q)}\right)^{(p,q)},$$
as can be verified by a straightforward calculation
that we omit, and we are done.
\qed

\subsubsection{Reparameterization invariance}
\label{subsubsection:ReparameterizationInvariance}
It is easily verified that for any $\tau\in k((t))$ of winding number
$1$, the operation 
$$(f(t)\mapsto f(\tau)):k((t))\rightarrow k((t))\;\;\;
(\mbox{``substitution of $\tau$ for $t$''})$$
is a $k$-linear automorphism of $k((t))$ belonging to the restricted
general linear group
$G^{k((t))}_{k[[t]]}$. Via the commutator interpretation 
provided by Theorem~\ref{Theorem:CommutatorInterpretation},
it follows that the Contou-Carr\`{e}re symbol is invariant under
reparameterization of
$k((t))$.

\subsubsection{Recovery of the tame symbol and the residue}
\label{subsubsection:ReductionRemark}
If $k$ is a field, the Contou-Carr\`{e}re symbol obviously reduces to
the tame symbol. It is possible also to recover the residue
from the Contou-Carr\`{e}re symbol, as follows. Take
$k=F[\epsilon]/(\epsilon^3)$ where
$F$ is any field. Then we have
$$\langle 1-\epsilon f,1-\epsilon g\rangle\equiv
1-\epsilon^2\Res_{t=0}(g\,df)\bmod{\epsilon^3}
$$
for all $f,g\in F((t))$ as can be verified by a straightforward
calculation. This last
observation suggests an interpretation of our work as the
``integrated version'' of Tate's Lie-theoretic theory
\cite{Tate}. We wonder how
Beilinson's multidimensional generalization \cite{Beilinson} of Tate's
theory
might analogously be integrated.

\subsubsection{The case in which $k$ is a
$\QQ$-algebra}\label{subsubsection:ExpLogEasy}
Suppose $k$ is a
$\QQ$-algebra. Let $f\in 1+m((t))$
and $g\in
k((t))^\times$ be given. 
We have
$$\langle f,g\rangle=\exp(\Res_{t=0}\log f\cdot d\log g)$$
as can be verified by a straightforward
calculation. This is quite similar in form to the
commutator formula 
given by Segal-Wilson \cite[Prop.\ 3.6]{SegalWilson}.

\section{Reciprocity laws on curves}
\subsection{The common setting for the reciprocity laws}

\subsubsection{Basic data}
Let $F$ be an algebraically
closed field. Let $X/F$ be a nonsingular complete algebraic curve.
Let $S$ be a finite nonempty set of closed points of $X$.
For any ring or
group
$A$, put
$$A^S:=\{(a_s)_{s\in S}\vert a_s\in A\}=
(\mbox{product of copies of $A$ indexed by $S$}).$$
For each $s\in S$ select a uniformizer $\pi_s$ at $s$.

\subsubsection{Construction of $R_0$}
For each meromorphic function $f$
on
$X$ put
$$f^{(s)}:=\sum_i a_it^i\in F((t))$$
where
$$f=\sum_i a_i\pi_s^i\;\;\;(a_i\in F)$$
is the Laurent expansion of $f$ in powers of $\pi_s$.
Put
$$R_0:=\{(f^{(s)})_{s\in S}\in F((t))^S\mid f\in H^0(X\setminus
S,\OO_X)\}.$$
The $F$-algebra $R_0$ is a copy of the affine coordinate ring of
$X\setminus S$.  We take for granted that
$$\dim_F F[[t]]^S\cap R_0=\dim_F H^0(X,\OO_X)=1<\infty$$
and
$$\dim_F \frac{F((t))^S}{R_0\oplus F[[t]]^S}=
\dim_F H^1(X,\OO_X)=\mbox{genus of $X$}<\infty.$$
As in the papers \cite{AdCK}, \cite{PablosRomo}, \cite{Tate}, 
it is these finiteness statements from algebraic geometry that lead
ineluctably to reciprocity laws.

\subsubsection{Extension of scalars from $F$ to $k$}
We assume now that the artinian local ring $k$ taken as base for
the theory of determinant groupoids is a finite
$F$-algebra. We put 
$$R:= R_0\otimes_F k=(\mbox{$k$-span of $R_0$})\subset
k((t))^S,$$
and
$$V:=k((t))^S,\;\;\;V_+:=k[[t]]^S.$$ Note that the
$k$-modules $V_+\cap R$ and
$V/(R+V_+)$ are free of finite rank. 
Note that $R^\times$ acting in natural $k$-linear fashion on $V$
is contained in the restricted general linear group $G_{V_+}^V$.

\subsection{A reciprocity law for the Contou-Carr\`{e}re symbol}
\begin{Theorem}
\label{Theorem:CCR}
In the setting above, for all $f,g\in R^\times$ we have
$$\prod_{s\in S}\langle f_s,g_s\rangle=1.$$
\end{Theorem}
\proof 
Clearly there exists some
$F$-subspace
$M_0\subset F((t))^S$ commensurable to $F[[t]]^S$ such that 
$F((t))^S=M_0\oplus R_0$. Put $M:=M_0\otimes_F k$.
Then we have
$$V=M\oplus R,\;\;\;M\sim V_+$$ 
and hence we
have
$$\{f,g\}_{V_+}^V=\{f,g\}_{M}^V=\{f,g\}_{M}^V\{f,g\}_R^V=1.$$ 
The first two equalities are justified by the basic properties
enumerated in \S\ref{subsubsection:BasicProperties} and
the last by
Proposition~\ref{Proposition:Switcheroo}.  
We also have
$$\{f,g\}_{V_+}^V=
\prod_{s\in S}
\prod_{s'\in S}\{f_s,g_{s'}\}_{k[[t]]}^{k((t))}
=
(-1)^{(\sum_s w(f_s))(\sum_s w(g_s))}\cdot\prod_{s\in S}\langle
f_s,g_s\rangle^{-1}.$$ The first equality is justified 
the basic properties enumerated in
\S\ref{subsubsection:BasicProperties} and
Propositions~\ref{Proposition:Flab}.
 The second equality is justified by
Proposition~\ref{Proposition:SignRule} and
Theorem~\ref{Theorem:CommutatorInterpretation}.  Finally, we have
$$\sum_{s\in S} w(f_s)=\sum_{s\in S} w(f_s\bmod{m})=0$$
because the second sum is the degree of a principal divisor on $X$. The
result follows.
\qed

\subsubsection{Coordinate-independence of the local
symbol
$\langle f_s,g_s\rangle$}
In
\S\ref{subsubsection:ReparameterizationInvariance} we explained  
that the Contou-Carr\`{e}re symbol is reparameterization invariant.
In the present context, it follows that the
value of the symbol
$\langle f_s,g_s\rangle$ is independent of the choice $\pi_s$ of
uniformizer at $s$, and thus is coordinate-independent.

\subsubsection{Recovery of Weil reciprocity}
Take $k=F$.
In this case the Contou-Carr\`{e}re symbol reduces to the tame
symbol, and hence Theorem~\ref{Theorem:CCR} reduces to
Weil reciprocity.

\subsubsection{Recovery of sum-of-residues-equals-zero}
Take $k=F[\epsilon]/(\epsilon^3)$.
In this case, as explained in \S\ref{subsubsection:ReductionRemark},
the residue can be recovered from the Contou-Carr\`{e}re symbol,
and hence
sum-of-residues-equals-zero
 can be recovered from Theorem~\ref{Theorem:CCR}.

\pagebreak
\subsection{Recovery of Witt's explicit reciprocity law}
\label{subsection:WittRecovery}

\subsubsection{Quick review of Witt vectors}
Witt vectors were introduced in Witt's paper 
\cite{Witt}. 
The basic theory---proofs omitted---takes the following form. 
For proofs, we recommend to the reader the exercises on this topic
in Lang's algebra text
\cite{Lang}. Let
$$\{\epsilon\}\coprod \{\xbold_i,\ybold_i\}_{i=1}^\infty$$
be a family of independent variables.
Write
$$\prod_{i=1}^\infty ((1-\xbold_i\epsilon^i)(1-\ybold_i\epsilon^i))
=\prod_{i=1}^\infty (1-\Abold_i\epsilon^i),\;\;\;
$$
$$
\prod_{i=1}^\infty \prod_{j=1}^\infty
\left(1-\xbold_i^{j/(i,j)}\ybold_j^{i/(i,j)}\epsilon^{ij/(i,j)}\right)^{(i,j)}=
\prod_{i=1}^\infty(1-\Mbold_i\epsilon^i),$$
thereby defining families of polynomials
$$\left\{\Abold_n,\Mbold_n\in \ZZ\left[\{\xbold_i,\ybold_i\}_{i\mid
n}\right]\right\}_{n=1}^\infty.$$
For any commutative ring $A$ with unit
and finite subset $\Delta$ of the set of positive integers closed under
passage to divisors, let
$\Witt_\Delta(A)$ denote the set of vectors with entries in $A$ indexed
by
$\Delta$.  It can be shown that the $\Abold$'s and $\Mbold$'s define
addition and multiplication laws with respect to which
$\Witt_\Delta(A)$ becomes a commutative ring with unit, functorially
in commutative rings $A$ with unit. Below we do not actually need
to use the multiplication law in $\Witt_\Delta(A)$ but we mention
its definition because of its close relationship to the 
definition of the Contou-Carr\`{e}re symbol.

\subsubsection{Ghost coordinates}
Consider the family of polynomials
$$\left\{\tilde{\xbold}_n:=\sum_{d\mid
n}d\xbold_{d}^{n/d}\right\}_{n=1}^\infty.$$
These polynomials are characterized by the power series identity
$$-\log \prod_{i=1}^\infty (1-\xbold_i\epsilon^i)
=\sum_{\nu=1}^\infty \frac{\tilde{\xbold}_\nu \epsilon^\nu}{\nu}.
$$
Given a ring
$A$, a finite subset
$\Delta$ of the set of positive integers closed under passage to
divisors,
$x=(x_i\in A)_{i\in \Delta}\in
\Witt_\Delta(A)$ and an integer $i\in \Delta$, we write $\tilde{x}_i$
for the result of substituting $x_d$ for $\xbold_d$ in
$\tilde{\xbold}_i$ for all $d\in \Delta$, 
and we call $\tilde{x}_i$ the
{\em ghost coordinate} of $x$ indexed by $i$; in this context,
for emphasis, we say that $x_i$ is the {\em live coordinate} of $x$
indexed by $i$.
Addition and multiplication have a
very simple expression in ghost
coordinates:
$$\tilde{\Abold}_n=\tilde{\xbold}_n+
\tilde{\ybold}_n,\;\;\;
\tilde{\Mbold}_n=\tilde{\xbold}_n\tilde{\ybold}_n.$$
Clearly each variable $\xbold_n$ has a unique expansion
as a polynomial in the $\tilde{\xbold}$'s with coefficients
in $\QQ$. It follows that for any $\QQ$-algebra $A$
the ring $\Witt_\Delta(A)$ decomposes in ghost
coordinates as a product of copies of
$A$ indexed by
$\Delta$.  But in general it is not
possible to write $\xbold_n$ as a polynomial
in the $\tilde{\xbold}$'s with integral coefficients,
and hence in general $\Witt_\Delta(A)$
depends in a complicated way on $A$. 

\subsubsection{Remark}
The focus in arithmetical applications of Witt vectors is usually on
the case $\Delta\subset\{1,p,p^2,\dots\}$ for some rational prime
$p$. In this case, for example, we have the striking
fact that
$$\Witt_{\{1,p,\dots,p^{n-1}\}}(\ZZ/p\ZZ)=\ZZ/p^n\ZZ.$$
But it is simpler to deal 
with the ring schemes of the form 
$$\Witt_{\leq
N}:=\Witt_{\{1,\dots,N\}}.$$
The additive group scheme underlying the ring scheme $\Witt_{\leq N}$
is fairly easy to handle  because the map
$$x=(x_i)_{i=1}^N\mapsto \prod_{i=1}^N
(1-x_i\epsilon)\bmod{\epsilon^{N+1}}$$ identifies
the additive group underlying $\Witt_{\leq N}(A)$
with the group of units in $A[\epsilon]/(\epsilon^{N+1})$
congruent to $1$ modulo $(\epsilon)$ functorially in commutative rings
$A$ with unit.
 Since the ring $\Witt_{\{1,p,\dots,p^{n-1}\}}(A)$ is a quotient
of the ring
$\Witt_{\leq p^{n-1}}(A)$ functorially in
$A$, it turns out that we really lose no generality by thus restricting
our focus.

\subsubsection{Definition of the symbol $\Res^{\Witt}_{\leq N}$}
Let us turn our attention back to the setting of
Theorem~\ref{Theorem:CCR}.
Fix a positive integer $N$. We
now take
$$k:=F[\epsilon]/(\epsilon^{N+1}).$$
 We define a pairing
$$\Res^\Witt_{\leq N}(\cdot,\cdot):F((t))^\times
\times \Witt_{\leq N}(F((t)))\rightarrow
\Witt_{\leq N}(F)$$
by the rule
$$\left\langle f,\prod_{i=1}^N\left(1-x_i\epsilon^i\right)\right\rangle
\equiv \prod_{i=1}^N
\left(1-\epsilon^i\left(\Res^\Witt_{\leq N}(f,x)\right)_i\right)
\bmod{(\epsilon^{N+1})}$$
where $\langle\cdot,\cdot\rangle$ is the Contou-Carr\`{e}re symbol.

\subsubsection{Comparison with Witt's original definition} The pairing
$\Res^{\Witt}_{\leq N}$ is essentially the pairing introduced in
Witt's paper \cite{Witt}. Without giving full details, we briefly
explain this point as follows. Assume for the moment that $F$ is of
characteristic zero so that we can talk about ghost coordinates, and
recall the remark of \S\ref{subsubsection:ExpLogEasy}. We have
$$\begin{array}{cl}
&\displaystyle\log \left\langle
f,\prod_{i=1}^N\left(1-\epsilon^ix_i\right)\right\rangle\\\\
\equiv&\displaystyle\Res_{t=0}
\left(-\log\left(\prod_{i=1}^{N}(1-\epsilon^i
x_i)\right)\cdot\frac{d\log f}{f}\right)\\\\
\equiv&\displaystyle\Res_{t=0}\left(
\sum_{i=1}^N\frac{\tilde{x}_i\epsilon^i}{i}\cdot\frac{d\log
f}{f}\right)\\\\
\equiv&\displaystyle\sum_{i=1}^N\left(\Res_{t=0}
\left(\tilde{x}_i\frac{df}{f}\right)\right)\frac{\epsilon^i}{i}.
\end{array}\bmod{(\epsilon^{N+1})}.
$$
In other words, we have
$$\widetilde{\Res^{\Witt}_{\leq N}(f,x)}_i=
\Res_{t=0}\left(\tilde{x}_i\frac{df}{f}\right)$$
for $i=1,\dots,N$.
This last formula
for $i=1,p,\dots,p^{n-1}$ and $F$ of characteristic $p>0$
is exactly
the expression in ghost coordinates 
of the rule used by Witt 
\cite[p. 130]{Witt} to define his pairing. But {\em a priori} Witt's
definition is no good in characteristic $p$ because live
coordinates cannot in general be expressed as polynomials in the ghost
coordinates with $p$-integral coefficients. Nevertheless,
Witt succeeds in
making the definition
rigorous  (see 
\cite[Satz 4, p.\ 130]{Witt}) by proving that the corresponding expression 
in live coordinates of his pairing is ``denominator-free'' and
hence does make sense in arbitrary characteristic.  
Of course, with the hindsight afforded by Theorem~\ref{Theorem:CCR},
Witt's Satz~4 is not very difficult to check.

\subsubsection{Reciprocity for the  symbol $\Res^\Witt_{\leq N}$}
We return to the situation in which the algebraically closed field $F$
may be of any characteristic. As above, we fix a positive integer $N$. 
We then have
$$\sum_{s\in S}\Res^{\Witt}_{\leq N}(f_s,x_s)=0$$
for all 
$$f\in R_0^\times\subset F((t))^{\times S},\;\;\;x\in
\Witt_{\leq N}(R_0)
\subset \Witt_{\leq N}(k((t))^S)=
\Witt_{\leq N}(k((t)))^S$$
by Theorem~\ref{Theorem:CCR} and the definitions.
We emphasize that the addition 
is to be
performed in the group $\Witt_{\leq N}(F)$. From this last formula 
in the case that $F$ is of characteristic $p$
and $N=p^{n-1}$, 
it is not difficult to deduce
the reciprocity law stated without proof on the last page of Witt's
paper
\cite{Witt} in the case of an algebraically closed ground
field of characteristic $p$. We omit further details.

\section{Acknowledgements}
The first-named author wishes to thank the mathematics department
at the University of Salamanca for its hospitality
during a two week period in May 2002 when much of the work 
on this paper was done.
The first-named author became acquainted with determinant groupoids 
from the perspective of mathematical physics through some
lectures on algebraic models of loop spaces by M.\ Kapranov 
in Spring 2002 at the University of Minnesota,
and gratefully acknowledges that influence.
The first-named author also wishes to thank the mathematics
department at the University of Arizona for the opportunity to
lecture on preliminary versions of some of these ideas at
the Winter School of March 2000.


\begin{thebibliography}{MM}



\bibitem{AdCK}
Arbarello, E.; De Concini, C.; Kac, V. G.:
{\em The infinite wedge representation and the reciprocity law 
for algebraic curves.} Theta
functions---Bowdoin 1987, Part 1 (Brunswick, ME, 1987), 171--190,
 Proc. Sympos. Pure Math., \textbf{49}, Part 1, Amer. Math. Soc., Providence, RI, 1989.











\bibitem{Beilinson}
Beilinson, A. A.:
{\em Residues and ad\`{e}les. (Russian)}
Funktsional. Anal. i Prilozhen. \textbf{14} (1980), no. 1, 44--45. 
\{English translation: Functional Anal.\ Appl.\ \textbf{14} 
(1980), no. 1, 34--35.\}


\bibitem{ContouCarrere} 
Contou-Carr\`{e}re, C.:
{\em Jacobienne locale, groupe de bivecteurs de Witt universel, et symbole 
mod\'er\'e. (French) [Local Jacobian,
universal Witt bivector group, and tame symbol]}
C.\ R.\ Acad.\ Sci.\ Paris S\'{e}r. I Math. \textbf{318} 
(1994), no. 8, 743--746. 



\bibitem{Lang}
Lang, S.: {\em Algebra. Revised third edition.} 
Graduate Texts in Mathematics, \textbf{211}. Springer-Verlag, New York, 2002. 




\bibitem{Matsumura}
Matsumura, H.: {\em Commutative algebra. Second edition.}
 Mathematics Lecture Note Series, \textbf{56}. 
Benjamin/Cummings Publishing
Co., Inc., Reading, Mass., 1980. 


\bibitem{PablosRomo} Pablos Romo, F.: {\em On the Tame Symbol of an
Algebraic Curve.} Comm.\ Algebra. \textbf{30}(2002), no.\ 9, 
4349--4368.





\bibitem{SegalWilson}
Segal, G.; Wilson, G.:
{\em Loop groups and equations of KdV type.}
 Inst.\ Hautes Études Sci.\ Publ.\ Math.\ \textbf{61}(1985), 5--65.

\bibitem{Tate}
Tate, J.: {\em Residues of differentials on curves.}
 Ann.\ Sci.\ Ecole Norm.\ Sup.\ (4) 1 1968 149--159.

\bibitem{Witt} Witt, E.:
{\em Zyklische K\"{o}rper und Algebren der Charakteristik $p$ vom
Grad $p^n$.} J.\ Reine  Angew.\ Math.\ {\bf 176}(1937), 126--140.


\end{thebibliography}
\end{document}